\documentclass[reqno,10pt]{amsart}
\usepackage{amsmath,amssymb}
\usepackage{latexsym}
\usepackage[dvips]{graphics}
\usepackage{epsfig}
\input amssym.def
\input amssym.tex

\newcommand{\zsum}[1]{ \sum_{#1 = 0}^{p-1} }
\newcommand{\osum}[1]{ \sum_{#1 = 1}^{p-1} }
\renewcommand{\i}{{\mathrm{i}}} % sqrt -1

\newcommand{\ep}[1]{ \mathbf{e}_p\left(#1\right) }

\newcommand\be{\begin{equation}}
\newcommand\ee{\end{equation}}
\newcommand\bea{\begin{eqnarray}}
\newcommand\eea{\end{eqnarray}}
\newcommand\bi{\begin{itemize}}
\newcommand\ei{\end{itemize}}
\newcommand\ben{\begin{enumerate}}
\newcommand\een{\end{enumerate}}
\newcommand\bc{\begin{center}}
\newcommand\ec{\end{center}}
\newcommand\ba{\begin{array}}
\newcommand\ea{\end{array}}
\newcommand{\R}{\ensuremath{\mathbb{R}}}

\newcommand{\Z}{\ensuremath{\mathbb{Z}}}
\newcommand{\Q}{\mathbb{Q}}

\newcommand{\foh}{\frac{1}{2}}  %onehalf
\newcommand{\js}[1]{{#1\overwithdelims()p}}

\newtheorem{thm}{Theorem}[section]

\newtheorem{lem}[thm]{Lemma}

\theoremstyle{definition}
\newtheorem{rek}[thm]{Remark}

\newcommand{\twocase}[5]{#1 \begin{cases} #2 & \text{{\rm #3}}\\ #4
&\text{\rm #5} \end{cases}   }

\newcommand{\oop}{\frac{1}{p}}

\begin{document}

\title[Variation in the number of points on
elliptic curves]{Variation in the number of points on elliptic
curves and applications to excess rank}

%\title{Michel's Second Moment Bound for Elliptic Curves is Sharp}

\author{Steven J. Miller}
\address{Department of Mathematics, Brown University, 151 Thayer
 Street,
Providence, RI 02912}
 \email{sjmiller@math.brown.edu}



\date{\today}

\begin{abstract}
Michel proved that for a one-parameter family of elliptic curves
over $\Q(T)$ with non-constant $j(T)$ that the second moment of
the number of solutions modulo $p$ is $p^2 + O(p^{3/2})$. We show
this bound is sharp by studying $y^2 = x^3 + Tx^2 + 1$. Lower
order terms for such moments in a family are related to lower
order terms in the $n$-level densities of Katz and Sarnak, which
describe the behavior of the zeros near the central point of the
associated $L$-functions. We conclude by investigating similar
families and show how the lower order terms in the second moment
may affect the expected bounds for the average rank of families in
numerical investigations.
\end{abstract}

\maketitle

\section{Introduction}

Let $\mathcal{E}$ be a one-parameter family of elliptic curves
(equivalently, an elliptic surface) over $\Q(T)$: \be y^2 +
a_1(T)xy + a_3(T)y \ = \ x^3 + a_2(T)x^2 + a_4(T)x + a_6(T), \ \ \
a_i(T) \in \Z[T]. \ee For each integer $t$ we have an elliptic
curve $E_t$ over $\Q$, with $N_t(p)$ the number of solutions
modulo $p$. Set $a_t(p) = p - N_t(p)$. If $a_1(T) = a_3(T) = 0$ we
have \be\label{eq:expansionatp} a_t(p) \ = \ -\sum_{x \bmod p}
\js{x^3 + a_2(t)x^2 + a_4(t)x + a_6(t)}. \ \ee We are interested
in evaluating the second moment for the family:
\begin{eqnarray}
A_{2,\mathcal{E}}(p) \ := \  \sum_{t \bmod p} a_t^2(p).
\end{eqnarray}
For one-parameter families of elliptic curves with $j(T)$
non-constant, Michel \cite{Mic} proves $A_{2,\mathcal{E}}(p) = p^2 +
O(p^{3/2})$ by using the Lefschetz-Groethendieck trace formula. We
show his result is sharp by constructing a family where the second
moment has a term of size $p^{3/2}$.

Moments of the number of solutions modulo $p$ provide enormous
amounts of information about the family. Rosen and Silverman
\cite{RoSi} prove a conjecture of Nagao (unconditionally for
rational elliptic surfaces, conditional on Tate's conjecture in
general) that the first moment is related to the rank of the
family over $\Q(T)$, denoted ${\rm rank}\ \mathcal{E}(\Q(T))$:
\begin{eqnarray} \lim_{X \to \infty} \frac{1}{X} \sum_{p \leq X}
\frac{\log p}p  \sum_{t \bmod p} a_t(p)  \ = \ -{\rm rank}\
\mathcal{E}(\Q(T)).
\end{eqnarray} In \cite{Mil1,ALM} it is shown how to construct
families with moderate rank by choosing the $a_i(T)$ so that the
first moment is computable and large.

Another application is in the connections between number theory
and random matrix theory \cite{KaSa1,KaSa2}. In showing the
behavior of the low lying zeros (zeros near the central point) of
$L$-functions of a family of elliptic curves agrees with that of
eigenvalues near $1$ of orthogonal groups, the only needed inputs
are the first and second moment of the number of solutions modulo
$p$; to first order all families of elliptic curves with rank $r$
over $\Q(T)$ agree with the same random matrix ensemble (see
\cite{Mil2}, though a similar result with a global rather than
local rescaling of the zeros is implicit in \cite{Sil2}). An
analysis of the lower order terms in the first and second moments
leads to breaking this universality; i.e., seeing lower-order
family dependent behavior in the low lying zeros. See
\cite{Mil1,Yo1} for more details on family dependent behavior. One
application of these correction terms is a refinement on
predicting the number of curves in a family with rank above the
family rank. While these corrections vanish in the limit of large
conductor, they lead to slight modifications of excess rank bounds
for conductors in the range accessible by computers.

In \S\ref{sec:2ndmomcalc} we determine the second moment of
$a_t(p)$ for a specific family, showing Michel's result is sharp.
In \S\ref{sec:excessrank} we analyze how the lower order terms in
Michel's theorem are related to bounds for the average rank of the
family.

\section{The Second Moment for
$y^2=x^3+Tx^2+1$}\label{sec:2ndmomcalc}

We may expand the Legendre symbol $\js{x}$ in
\eqref{eq:expansionatp} by
\begin{eqnarray}\label{eq:jsgaussum}
\js{x} \ = \  \frac1{G_p} \sum_{c=1}^{p-1} \js{c} \ep{cx}, \ \ \
\ep{x} \ = \ e^{2\pi \i x /p}.
\end{eqnarray}
Here $G_p = \sum_{a \bmod p} \js{a} \ep{a}$ is the Gauss sum,
which equals $\sqrt{p}$ for $p \equiv 1 \bmod 4$ and $\i\sqrt{p}$
for $p \equiv 3 \bmod 4$. See, for example, \cite{BEW}. Using
Gauss sums to evaluate Legendre sums is a common technique; we
sketch an alternate approach which avoids Gauss sums in Remark
\ref{rek:cooknog}.

For the family $\mathcal{E}: y^2 = x^3 + Tx^2 + 1$, $j(T) = -
\frac{256T^6}{4T^3+27}$ and thus by Michel's Theorem
$A_{2,\mathcal{E}}(p) = p^2 + O(p^{3/2})$. We determine an exact
formula for $A_{2,\mathcal{E}}(p)$:

\begin{thm} For the one-parameter family $\mathcal{E}: y^2 = x^3 + Tx^2 + 1$
over $\Q(T)$, for $p>2$ the second moment of $a_t(p)$ is \be
A_{2,\mathcal{E}}(p) \ = \ \sum_{t\bmod p} a_t(p)^2 \ = \  p^2 -
n_{3,2,p}p - 1 + p\sum_{x \bmod p} \js{4x^3+1}, \ee where
$n_{3,2,p}$ denotes the number of cube roots of $2$ modulo $p$. For
any $[a,b] \subset [-2,2]$ there are infinitely many primes $p
\equiv 1 \bmod 3$ such that \be A_{2,\mathcal{E}}(p) - \left(p^2-
n_{3,2,p}p - 1\right)\ \in\ [a\cdot p^{3/2},b\cdot p^{3/2}].\ee
\end{thm}

\begin{proof} Combining \eqref{eq:expansionatp} and \eqref{eq:jsgaussum}
yields
\begin{align}\label{eq:expfamxycd}
& A_{2,\mathcal{E}}(p) \ = \ \sum_{t \bmod p}\ \sum_{x \bmod p}\
\sum_{y \bmod p} \js{x^3+1+x^2t} \js{y^3+1+y^2t} \nonumber\\ & \ \
\ =\ \sum_{x,y \bmod p}\ \sum_{c,d = 1}^{p-1} \frac{1}{p} \js{cd}
\ep{c(x^3+1)-d(y^3+1)} \sum_{t \bmod p} \ep{(cx^2-dy^2)t};
\end{align} above we used the complex conjugate of
\eqref{eq:jsgaussum} in expanding $\js{y^3+1+y^2t}$. The two Gauss
sum expansions give $\frac1{G_p\overline{G_p}} = \frac1p$. It will
be convenient to set $g(x,y) = (x-y)(x^2y^2 - (x+y))$.

Note $c$ and $d$ are invertible modulo $p$ in \eqref{eq:expfamxycd}.
If the numerator in the $t$-exponential is non-zero, the $t$-sum
vanishes. Thus it suffices to study \eqref{eq:expfamxycd} for $cx^2
\equiv dy^2 \bmod p$. If exactly one of $x$ and $y$ vanishes, then
$cx^2 \not\equiv dy^2 \bmod p$. Hence $x$ and $y$ are both zero or
non-zero. If both are zero the $t$-sum gives $p$, the $c$-sum gives
$G_p$, the $d$-sum gives $\overline{G}_p$, for a total contribution
of $p$. If $x$ and $y$ are non-zero then we must have $d \equiv
cx^2y^{-2} \bmod p$. The $t$-sum gives $p$. Thus
\eqref{eq:expfamxycd} is
\begin{eqnarray}\label{eq:startA2Ep}
A_{2,\mathcal{E}}(p) &\ = \ & p+ \sum_{x,y =1}^{p-1} \sum_{c =
1}^{p-1} \frac{1}{p} \js{x^2y^2} \ep{cy^{-2}(x^3y^2 + y^2 - x^2y^3
- x^2)}p \nonumber\\ &=& p+ \sum_{x,y =1}^{p-1} \sum_{c = 1}^{p-1}
\ep{cy^{-2}(x-y)(x^2y^2 - (x+y))} \nonumber\\ &=& p+ \sum_{x,y
=1}^{p-1} \sum_{c = 0}^{p-1}
\ep{cy^{-2}g(x,y)} - \sum_{x,y =1}^{p-1} 1 \nonumber\\
&=& p+ \sum_{x,y =1}^{p-1} \sum_{c = 0}^{p-1} \ep{cy^{-2}g(x,y)} -
(p-1)^2 \nonumber\\ & = & p\cdot\#\{ x, y \not\equiv 0 \bmod p:
g(x,y) \equiv 0 \bmod p\} + p - (p-1)^2,
\end{eqnarray} where the last equality follows from the fact that if
$g(x,y) \equiv 0 \bmod p$ then the $c$-sum is $p$ and otherwise it
is $0$. We are left with counting how often $g(x,y) \equiv 0\bmod
p$ for $x$, $y$ non-zero.

Whenever $x = y$ then $g(x,y) \equiv 0 \bmod p$; therefore there
are $p-1$ solutions from $x = y$. Consider now $x^2y^2 \equiv x+y
\bmod p$, which we may rewrite as a quadratic in $y$: $x^2 y^2 - y
- x$ $ \equiv 0\bmod p$. By the Quadratic Formula modulo $p$
(recall $p$ is odd), if the discriminant $4x^3+1$ is a non-zero
square modulo $p$ there are two distinct roots, if it is not a
square modulo $p$ there are no roots, and if the discriminant
vanishes there is one root. Equivalently, if $\js{4x^3 + 1} = 1$
(respectively, $-1$ or $0$) there are two (respectively, none or
one) solutions to $x^2y^2 \equiv x+y\bmod p$.

Recall neither $x$ nor $y$ is allowed to be zero. If $y = 0$ then
$x^2y^2 \equiv x+y\bmod p$ reduces to $x=0$. Hence our solutions
have $x, y \not\equiv 0\bmod p$, and for a non-zero $x$ the number
of non-zero $y$ with $x^2y^2 \equiv x+y \bmod p$ is $1 +
\js{4x^3+1}$. Hence the number of non-zero pairs with $x^2y^2
\equiv x+y\bmod p$ is
\begin{eqnarray}
\sum_{x=1}^{p-1} \left(1 + \js{4x^3+1} \right)\ =\ p-1 + \sum_{x
\bmod p} \js{4x^3 + 1} - 1.
\end{eqnarray}

We must be careful about double counting solutions. If both $x-y
\equiv 0 \bmod p$ and $x^2y^2 \equiv x+y \bmod p$, then we find
$x^4 \equiv 2x \bmod p$. As $x \not\equiv 0 \bmod p$, we obtain
$x^3 \equiv 2 \bmod p$. We have double counted all pairs $(x,x)$
with $x$ a cube root of $2$ modulo $p$. Let $n_{3,2,p}$ denote the
number of cube roots of $2$ modulo $p$; $|n_{3,2,p}| \le 3$. We
have shown
\begin{align} & \#\{ x, y \in \{1,\dots,p-1\}: g(x,y) \equiv 0
\bmod p\} \nonumber\\ & \ \ \ \ \ \ \ = \  (p-1) + \left(p-1 +
\sum_{x \bmod p} \js{4x^3 + 1} - 1\right) - n_{3,2,p} \nonumber\\
& \ \ \ \ \ \ \ = \ 2p - 3 - n_{3,2,p} + \sum_{x \bmod p} \js{4x^3
+ 1}.
\end{align}
Thus \eqref{eq:startA2Ep} becomes
\begin{eqnarray}
A_{2,\mathcal{E}}(p) & = & p\left(2p - 3 - n_{3,2,p} + \sum_{x \bmod
p} \js{4x^3 + 1}\right) + p - (p-1)^2  \nonumber\\ &=& p^2 -
n_{3,2,p}p - 1 + p\sum_{x \bmod p} \js{4x^3+1}.
\end{eqnarray}

To complete the analysis, we need to determine the size of $\sum_{x
\bmod p} \js{4x^3+1}$. Note this is the number of solutions modulo
$p$ to the elliptic curve $y^2 = 4x^3+1$, and this curve is
equivalent to $E: y^2 = x^3 + 16$. This curve has analytic rank $0$,
as can be seen from $L(E,1) \sim .5968$. It has complex
multiplication, and for $p \equiv 2 \bmod 3$, $a_E(p) = 0$. Write
$a_E(p) = 2\sqrt{p}\cos \theta_{E,p}$. As $E$ has complex
multiplication, for $p\equiv 1 \bmod 3$ the distribution of the
angles $\theta_{E,p}$ is known; all we need is that for any
$[\Theta,\Theta'] \subset [0,\pi]$, a positive percent of the time
$\theta_{E,p} \in [\Theta,\Theta']$. This implies that a typical
$a_E(p)$ is of size $\sqrt{p}$ if $p\equiv 1 \bmod 3$, and hence for
any $[a,b] \subset [-2,2]$ we can find infinitely many primes $p$
with\be A_{2,\mathcal{E}}(p) - \left(p^2- n_{3,2,p}p - 1\right)\
\in\ [a\cdot p^{3/2},b\cdot p^{3/2}].\ee
\end{proof}

Note that if $p \equiv 2 \bmod 3$, as $x \mapsto x^3 \bmod p$ is an
automorphism then $n_{3,2,p} = 1$ and $a_E(p) = 0$. Thus, at least
half the time, $A_{2,\mathcal{E}}(p) = p^2 - p - 1$.

\begin{rek}\label{rek:cooknog}
A few words should be said about how we cooked up this family. If
instead of $y^2 = x^3 + T x^2 + 1$ we had $y^2 = x^3 + T x + 1$, we
would have found the condition $d \equiv cxy^{-1} \bmod p$. As we
have $\js{cd}$ this would lead to $\js{c^2}\js{xy}$ times a similar
$c$-exponential. It would not suffice to determine how often a
similar $g(x,y)$ vanished; we would need to know the value of
$\js{xy}$. Our analysis was greatly aided by the presence of
$\js{x^2 y^2}$. We also want to change the order of summation and do
the $t$-sum first, which basically forces our family to be at most
quadratic in $t$, and such that $g(x,y)$ factors easily. Instead of
expanding by using Gauss sums, we could write the product of
Legendre symbols \eqref{eq:expfamxycd} as the Legendre symbol of
$h(x,y,t)$, where $h$ is quadratic in $t$ with leading term $x^2y^2
t^2$: \begin{align} & \js{x^3+1+x^2t} \js{y^3+1+y^2t} \nonumber\\ &
\ \ \ \ \ = \ \js{x^2 y^2 \cdot t^2 +
\left(y^2(x^3+1)+x^2(y^3+1)\right)\cdot t - (x^3+1)(y^3+1)}.
\end{align} We execute the $t$-sum first. Quadratic Legendre sums
are easily determined; what matters is the discriminant modulo
$p$. After some algebra we find the discriminant is $g(x,y)^2$
(with $g(x,y)$ as before), and then the argument proceeds
identically. See \cite{Mil1,ALM} for more on determining tractable
families where the summation can be done in closed form. These
families will be quadratic in $t$, although not necessarily in
Weierstrass form.
\end{rek}

\section{Other Families and Applications to Excess Rank}\label{sec:excessrank}

We give some additional examples of families where the first and
second moments can be determined exactly; see \cite{Mil1} for the
calculations (though we provide calculations of a representative
set of these families in Appendices \ref{sec:fam1} through
\ref{sec:fam3}).

Recall $n_{3,2,p}$ denotes the number of cube roots of $2$ modulo
$p$, and set $c_0(p) = \js{-3} + \js{3}$, $c_1(p) = \left[\sum_{x
\bmod p} \js{x^3-x}\right]^2$ and $c_{3/2}(p) = p\sum_{x \bmod p}
\js{4x^3+1}$.

\begin{tabular}{lcl} \\
${\ \ \ \ \ \ \ \mbox{Family}  \ \ \ \ \ \ \ }$ & ${\ \ \ \ \ \ \
A_{1,\mathcal{E}}(p) \ \ \ \ \ \ \ }$ & ${\ \ \ \ \ \ \
A_{2,\mathcal{E}}(p) \ \ \ \ \ \ \ }$
\\ \hline

$y^2 = x^3 + Sx + T$ & $\ \ 0$    & $p^3 - p^2$ \\

$y^2 = x^3 + 2^4(-3)^3(9T+1)^2$ & $\ \ 0$ & $\Big\{ {2p^2 - 2p \ \
\ p
\equiv 2 \bmod 3\atop \ 0 \ \ \ \ \ \ \ \ \ p \equiv 1 \bmod 3 }$ \\

$y^2 = x^3 \pm 4(4T+2)x$ & $\ \ 0$ & $\Big\{ {2p^2 - 2p \ \ \  p
\equiv
1 \bmod 3 \atop \ 0 \ \ \ \ \ \ \ \ \ p \equiv 3 \bmod 3}$ \\

$y^2 = x^3 + (T+1)x^2 +Tx$ & $\ \ 0$ & $p^2 - 2p - 1$ \\

$y^2 = x^3 + x^2 + 2T+1$ & $\ \ 0$ & $p^2 - 2p - \js{-3}$ \\

$y^2 = x^3 + Tx^2 + 1$ & $-p$ & $p^2 -
n_{3,2,p}p - 1 + c_{3/2}(p)$ \\

$y^2 = x^3 - T^2x + T^2$ & $-2p$ & $p^2 - p - c_1(p) - c_0(p)$ \\

$y^2 = x^3 - T^2x + T^4$ & $-2p$ & $p^2 - p - c_1(p) - c_0(p)$

\end{tabular}\\

The first family is the family of all elliptic curves; it is a two
parameter family and we expect the main term of its second moment to
be $p^3$. Note that except for our family $y^2 = x^3 + Tx^2 + 1$,
all the families $\mathcal{E}$ have $A_{2,\mathcal{E}}(p) = p^2 -
h(p)p + O(1)$, where $h(p)$ is non-negative. Further, many of the
families have $h(p) = m_\mathcal{E} > 0$. Note $c_1(p)$ is the
square of the coefficients from an elliptic curve with complex
multiplication. It is non-negative and of size $p$ for $p\not\equiv
3\bmod 4$, and zero for $p\equiv 1 \bmod 4$ (send $x \mapsto -x
\bmod p$ and note $\js{-1}=-1$). It is somewhat remarkable that all
these families have a correction to the main term in Michel's
theorem in the same direction, and we analyze the consequence this
has on the average rank. For our family which has a $p^{3/2}$ term,
note that on average this term is zero and the $p$ term is negative.

Consider a one-parameter family of elliptic curves $\mathcal{E}$
of rank $r$ over $\Q(T)$. With our normalizations, under GRH the
non-trivial zeros of $E_t$ are $1+i\gamma_t$, $\gamma_t\in\R$. We
typically study $t\in [N,2N]$ with $N\to\infty$. Let $C_t$ be the
conductor of the elliptic curve $E_t$, and let $\log R = \frac1N
\sum_{t=N}^{2N} \log C_t$ be the average log-conductor. For many
families there is an integer $a$ such that $\log C_t \sim \log
N^a$ for most curves; this is true for the families listed above
(see \cite{Mil1} for the calculations). Assuming the Birch and
Swinnerton-Dyer conjecture, by Silverman's specialization theorem
eventually all curves $E_t$ have rank at least $r$, and under
natural standard conjectures (see \cite{He}) a typical family will
have equidistribution of signs of the functional equations. What
is typically seen in studying the ranks of curves in a family is
that roughly $30\%$ have rank $r$ and $20\%$ rank $r+2$, while
about $48\%$ have rank $r+1$ and $2\%$ rank $r+3$. Random matrix
theory predicts that in the limit $50\%$ should be rank $r$ and
$50\%$ rank $r+1$ for an average rank of $r+\foh$, markedly
different from the observed (approximately) $r + \foh + .40$. See
\cite{Fe1,Fe2,Wa} for numerical investigations and
\cite{Br,H-B,FP,Mic,Sil2,Yo2} for theoretical bounds of the
average rank.

The excess rank question is whether this disagreement persists or is
a result of small data. We often expect the rate of convergence for
problems such as this to be like the logarithm of the conductors. As
the conductors are often at most $10^{12}$, it is reasonable to
believe the data is misleading (especially as random matrix theory
predicts sub-families of higher rank of size $N^{3/4}$, and for
small $N$ such families are a noticeable percentage; see for example
\cite{CKRS,DFK,Go,GM,Mai,Ono,RoSi,ST,Yu} for discussions of random
matrix predictions and results from number theory).

For an even Schwartz test function $\phi$ with ${\rm
supp}(\widehat{\phi}) \subset (-\sigma,\sigma)$, the $1$-level
density (which is basically just the sum of the explicit formula for
each curve) is defined by
\begin{align}\label{eqoneleveldensityexpansion}
& \frac1N \sum_{t=N}^{2N}\sum_{\gamma_t}
\phi\left(\gamma_t\frac{\log R}{2\pi}\right) \ = \
\widehat{\phi}(0) + \phi(0) - \ \frac2N \sum_{t=N}^{2N} \sum_p
\frac{\log p}{\log R} \oop \widehat{\phi}\left(\frac{\log p}{\log
R} \right)
a_t(p) \nonumber\\
& \ \ \ - \ \frac2N \sum_{t=N}^{2N} \sum_p \frac{\log p}{\log R}
\frac{1}{p^2} \widehat{\phi}\left(\frac{2\log p}{\log R} \right)
a_t(p)^2 + O\left(\frac{\log\log R}{\log R}\right).
\end{align} If $\phi$ is non-negative, we obtain a bound for the
average rank in the family by restricting the sum to be only over
zeros at the central point. The error $O\left(\frac{\log \log
R}{\log R}\right)$ comes from trivial estimation and ignores
probable cancellation, and we expect $O\left(\frac1{\log
R}\right)$ or smaller to be the correct magnitude. For most
families $\log R \sim \log N^a$ for some integer $a$.

The main term of the first and second moments of the $a_t(p)$ give
$r\phi(0)$ and $-\foh \phi(0)$, respectively, in
\eqref{eqoneleveldensityexpansion}. Assume the second moment of
$a_t(p)^2$ is $p^2 - m_\mathcal{E}p + O(1)$, $m_\mathcal{E} > 0$.
We have already handled the contribution from $p^2$, and $-
m_\mathcal{E}p$ contributes
\begin{eqnarray}\label{eq:psumforS2}
S_2 &\ \sim\ & \frac{-2}{N} \sum_p \frac{\log p}{\log
R}\widehat{\phi}\left(2\frac{\log p}{\log R}\right) \frac{1}{p^2}
\frac{N}{p}(-m_{\mathcal{E}}p) \nonumber\\ &=&
\frac{2m_{\mathcal{E}}}{\log R} \sum_p
\widehat{\phi}\left(2\frac{\log p}{\log R} \right) \frac{\log
p}{p^2}.
\end{eqnarray}
Thus there is a contribution of size $\frac{1}{\log R}$. A good
choice of test functions (see Appendix A of \cite{ILS}) is the
Fourier pair
\begin{eqnarray}
\phi(x) \ = \  \frac{\sin^2(2\pi \frac{\sigma}2 x)}{(2\pi x)^2}, \
\ \ \ \twocase{\widehat{\phi}(u) \ =\ }{\frac{\sigma - |u|}{4}}{if
$|u| \leq \sigma$}{0}{otherwise.}
\end{eqnarray}
Note $\phi(0) = \frac{\sigma^2}{4}$, $\widehat{\phi}(0) =
\frac{\sigma}4 = \frac{\phi(0)}{\sigma}$, and evaluating the prime
sum in \eqref{eq:psumforS2} gives
\begin{eqnarray}
S_2 & \ \sim \ & \left(\frac{.986}{\sigma} - \frac{2.966}{\sigma^2
\log R} \right)\frac{m_{\mathcal{E}}}{\log R}\ \phi(0).
\end{eqnarray}

Let $r_t$ denote the number of zeros of $E_t$ at the central point
(i.e., the analytic rank). Then up to our $O\left(\frac{\log\log
R}{\log R}\right)$ errors (which we think should be smaller), we
have
\begin{eqnarray}
\frac1N \sum_{t=N}^{2N} r_t\phi(0) & \leq & \frac{\phi(0)}{\sigma}
+ \left(r+\foh\right)\phi(0) + \left( \frac{.986}{\sigma} -
\frac{2.966}{\sigma^2
\log R} \right)\frac{m_{\mathcal{E}}}{\log R}\phi(0) \nonumber\\
\mbox{Ave Rank}_{[N,2N]}(\mathcal{E}) & \leq & \frac{1}{\sigma} +
r + \foh + \left( \frac{.986}{\sigma} - \frac{2.966}{\sigma^2 \log
R} \right)\frac{m_{\mathcal{E}}}{\log R}.
\end{eqnarray}

\begin{rek} The Density Conjecture states that the
$1$-level density (in the limit) should hold for all $\sigma$. In
that case, the lower order terms from the second moment will
\emph{not} contribute to the bound for the average rank, as their
contribution vanishes as $\sigma \to \infty$. Of course, the
agreement with random matrix theory is a statement about the limit
as $N\to\infty$; the correct finite-conductor model is still
unknown
\end{rek}

Let us examine the boost the $-m_\mathcal{E}p$ term from the
second moment gives to the upper bound for the average rank. As
remarked, if our $1$-level density were true for all $\sigma$
then there would be no contribution from the correction term to
the second sum, nor would the $\frac{1}{\sigma}$ term contribute,
and we would obtain the average rank is bounded by $r+\foh$.

Let us assume we know the $1$-level density up to $\sigma = 1$.
(This is well beyond the range of current technology; the best
result to date is for the family of all elliptic curves, where Young
\cite{Yo2} proves we may take any $\sigma < \frac79$). Assume
$m_{\mathcal{E}} = 1$. The $\frac{1}{\sigma}$ term would contribute
$1$, the lower correction would contribute $.03$ for conductors of
size $10^{12}$, and thus the average rank is bounded by $1 + r +
\foh + .03$ $= r + \foh + 1.03$. This is significantly higher than
Fermigier's observed $r + \foh + .40$.

If we were able to prove our $1$-level density for $\sigma = 2$,
then the $\frac{1}{\sigma}$ term would contribute $\foh$, and the
lower order correction would contribute $.02$ for conductors of size
$10^{12}$. Thus the average rank would be bounded by $\foh + r +
\foh + .02$ $= r + \foh + .52$. While the main error contribution is
from $\frac{1}{\sigma}$, there is still a noticeable effect from the
lower order terms in $A_{2,\mathcal{E}}(p)$. Moreover, we are now in
the ballpark of Fermigier's bound; of course, we were already there
without the potential correction term.

It seems hopeless to think about obtaining a $1$-level density for
any family of elliptic curves with support $\sigma =2$ or more.
Iwawniec, Luo and Sarnak \cite{ILS} obtain such large support for
families of weight $k$ cuspidal newforms of square-free level $N$,
but only because of great averaging formulas (the Bessel-Kloosterman
expansion in the Petersson formula) available for the family; the
corresponding averaging formulas for elliptic curves are much
weaker. We use the periodicity of $a_t(p)$ as a function of $t \bmod
p$ to analyze complete sums of the moments for each prime; however,
the error from the incomplete sum is bounded by Hasse and
contributes a large error (this problem is avoided in \cite{ILS}
because the Petersson formula gives us sums over a basis of
newforms, and there is no incomplete piece to be approximated). The
random matrix models for the behavior of the zeros near the central
point have been shown to hold as the conductors tend to infinity; in
the small conductor ranges investigated, it is not surprising that
there is disagreement. While it would be desirable to find a good
model for small conductors (similar to Keating and Snaith's
\cite{KeSn1,KeSn2} modeling zeros of $\zeta(s)$ at height $T$ by
$N\times N$ matrices with $N = \frac{\log T}{2\pi}$), we can
identify potential family dependent lower order terms in the
$1$-level density arising from lower order terms in the second
moment. For finite conductors these do lead to slightly larger
predicted upper bounds for the average rank in a family.

\appendix

\newpage
%\ \\ \ \\ \ \\ \ \\ \ \\ \ \\

\LARGE \textsf{The following appendices contain calculations for
the second moment of a representative set of the other families
mentioned in \S\ref{sec:excessrank}, and for completeness proofs
of standard quadratic Legendre sums. For general one-parameter
families of elliptic curves, there are not closed form expressions
for the moments. Our hope is that these families may be useful for
other investigations.}\normalsize

\ \\

\section{The Family $y^2 = x^3 + T^2$}\label{sec:fam1}

\begin{thm} For the family $\mathcal{E}: y^2 = x^3 + T^2$,
\be \twocase{A_{2,\mathcal{E}}(p) \ = \ }{2p^2-2p}{if $p \equiv 1
\bmod 3$}{0}{if $p \equiv 2 \bmod 3$.} \ee \end{thm}

\begin{proof} If $p \equiv 2 \bmod 3$ then $a_t^2(p) = 0$ as $x\mapsto x^3
\bmod p$ is an automorphism. Assume $p \equiv 1 \bmod 3$.
\begin{eqnarray}
A_{2,\mathcal{E}}(p) &\ = \ &  \sum_{t \bmod p}\ \sum_{x \bmod p}\
\sum_{y \bmod p} \js{x^3+t^2}\js{y^3+t^2} \nonumber\\ &=&
\sum_{t=1}^{p-1}
\sum_{x \bmod p}\ \sum_{y \bmod p} \js{x^3+t^2}\js{y^3+t^2} \nonumber\\
&=& \sum_{t=1}^{p-1}\ \sum_{x \bmod p}\ \sum_{y \bmod p} \js{t^4}
\js{tx^3+1}\js{ty^3+1} \nonumber\\ &=& \sum_{x \bmod p}\ \sum_{y
\bmod p}\ \sum_{t \bmod p} \js{tx^3+1}\js{ty^3+1} - p^2.
\end{eqnarray}
We use inclusion / exclusion to reduce to $xy \neq 0$. If $x=0$,
the $t$-sum vanishes unless $y=0$, in which case we get $p$.
Similarly if $y=0$, the $t$ and $x$-sums give $p$. We subtract the
doubly counted contribution from $x=y=0$, which gives $p$. Thus
\begin{eqnarray}
A_{2,\mathcal{E}}(p) &\ = \ & \sum_{x=1}^{p-1}\sum_{y=1}^{p-1}
\sum_{t (p)} \js{tx^3+1}\js{ty^3+1} + 2p -p - p^2.
\end{eqnarray}

By Lemma \ref{lablegsum}, the $t$-sum is $(p-1)\js{x^3y^3}$ if
$p|(x^3-y^3)^2$ and $-\js{x^3y^3}$ otherwise. As $p = 6m+1$, let
$g$ be a generator of the multiplicative group $\Z/p\Z$. Solving
$g^{3a} \equiv g^{3b}\bmod p$ yields $b = a$, $a + 2m$, or $a +
4m$. Thus, $x^3 \equiv y^3 \bmod p$ three times, and in each
instance $y$ equals $x$ times a square ($1$, $g^{2m}$, $g^{4m}$).
\begin{eqnarray}
A_{2,\mathcal{E}}(p) &=& \sum_{x=1}^{p-1} \sum_{y=1 \atop y^3
\equiv x^3}^{p-1} p - \sum_{x=1}^{p-1}\sum_{y=1}^{p-1} \js{x^3y^3}
+ p - p^2 \nonumber\\ &=& (p-1)3p + p - p^2 \nonumber\\ &=& 2p^2 -
2p.
\end{eqnarray} \end{proof}

\section{The Family $y^2 = x^3+x^2+T$}\label{sec:fam2}

\begin{thm} For the family $\mathcal{E}: y^2 = x^3 + x^2 + T$ we
have
\begin{eqnarray}
A_{2,\mathcal{E}}(p) \ = \ p^2 - 2p - p\js{-3}.
\end{eqnarray} \end{thm}

\begin{proof}
\begin{eqnarray}
A_{2,\mathcal{E}}(p) &\ = \ & \sum_{t=0}^{p-1}
\sum_{x=0}^{p-1}\sum_{y=0}^{p-1} \js{t + (x^3+x^2)}\js{t +
(y^3+y^2)} \nonumber\\ & = & \sum_{t=0}^{p-1}
\sum_{x=0}^{p-1}\sum_{y=0}^{p-1} \js{t^2 +
\left((x^3+x^2)+(y^3+y^2)\right)t +
(x^3+x^2)(y^3+y^2)}.\nonumber\\
\end{eqnarray} Let $\delta(x,y) = (x^3 + x^2) - (y^3+y^2)$; note
$\delta(x,y)^2$ is the discriminant of the quadratic regarded as a
function of $t$. The $t$-sum is $p-1$ if $p|\delta(x,y)$ and $-1$
otherwise. Note \be \delta(x,y) \ = \ (x-y)(y^2 + (x+1)y +
(x^2+x)). \ee The first factor is congruent to zero when $x = y$;
for fixed $x$, the discriminant of the second factor is $(x+1)^2 -
4(x^2+x) = 1 - 2x - 3x^2$. Thus the number of solutions of the
second factor, for fixed $x$, is $1 + \js{1-2x-3x^2}$. As the
discriminant of $1-2x-3x^2$ is $16$, summing over $x$ for $p > 2$
yields $p - \js{-3}$ by Lemma \ref{labquadlegsum}.

We must be careful about double counting. If both factors are
congruent to zero, then $3x^2 + 2x \equiv 0$, or $x \equiv 0,
-2\cdot 3^{-1}$. Hence we always double count two solutions.
\begin{eqnarray}
A_{2,\mathcal{E}}(p) &\ = \ & \Bigg[p + p -\js{-3} - 2\Bigg]p -
\sum_{x=0}^{p-1}\sum_{y=0}^{p-1} 1 \nonumber\\ &=& p^2 - 2p -
p\js{-3}.
\end{eqnarray}\end{proof}

\section{$y^2 = x^3 - T^2x + T^2$}\label{sec:fam3}

Consider the family $\mathcal{E}: y^2 = x^3 -T^2x + T^2$. We
calculate the first moment of $a_t(p)$, which shows the family is
of rank $2$ over $\Q(T)$, and then determine the second moment.

\begin{thm} For $\mathcal{E}: y^2 = x^3 -T^2x + T^2$,
$A_{1,\mathcal{E}}(p) = -2p$. Thus by Rosen and Silverman the
family has rank $2$ over $\Q(T)$. \end{thm}
\begin{proof}
\begin{eqnarray}
-A_{1,\mathcal{E}}(p) &=& -\sum_{t (p)} a_t(p) \ = \ \sum_{t (p)}
\sum_{x (p)} \js{x^3 -t^2x + t^2} \nonumber\\ &=& \sum_{t=1}^{p-1}
\sum_{x (p)} \js{x^3 -t^2x + t^2} \ = \ \sum_{t=1}^{p-1} \sum_{x
(p)} \js{t^3x^3 -t^3x + t^2} \nonumber\\ &=& \sum_{t=1}^{p-1}
\sum_{x (p)} \js{t^2} \js{t(x^3 -x) + 1} \nonumber\\ &=& \sum_{t
(p)} \sum_{x (p)} \js{t(x^3 - x) + 1} - \sum_{x (p)} \js{1}
\nonumber\\ &=& \sum_{t (p)} \sum_{x=0,\pm 1} \js{t(x^3 - x) + 1}
+ \sum_{t (p)} \sum_{x (p) \atop x \neq 0,\pm 1} \js{t(x^3 - x) +
1} - p \nonumber\\ &=& \sum_{t (p)} \sum_{x=0,\pm 1} \js{1} +
\sum_{x (p) \atop x \neq 0,\pm 1} \sum_{t (p)} \js{t + 1} - p
\nonumber\\ &=& 3p + 0 - p \ = \ 2p.
\end{eqnarray} \end{proof}

\begin{thm}For $\mathcal{E}: y^2 = x^3 -T^2x + T^2$,
\begin{eqnarray}
A_{2,\mathcal{E}}(p) = p^2 - p - \Bigg[\sum_{x (p)}
\js{(x^3-x)}\Bigg]^2 - \js{-3} - \js{3} = p^2 + O(p). \nonumber\\
\end{eqnarray}
\end{thm}

\begin{proof}
\begin{eqnarray}
A_{2,\mathcal{E}}(p) &=& \sum_{t (p)} a_t^2(p) \nonumber\\ &=&
\sum_{t (p)} \sum_{x,y (p)} \js{x^3 -t^2x + t^2} \js{y^3-t^2y+t^2}
\nonumber\\ &=& \sum_{t=1}^{p-1} \sum_{x,y (p)} \js{x^3 -t^2x +
t^2} \js{y^3-t^2y+t^2} \nonumber\\ &= & \sum_{t=1}^{p-1} \sum_{x,y
(p)} \js{t^3x^3 -t^3x + t^2} \js{t^3y^3-t^3y+t^2} \nonumber\\ &= &
\sum_{t=1}^{p-1} \sum_{x,y (p)} \js{t^4} \js{t(x^3-x)+1}
\js{t(y^3-y)+1} \nonumber\\  &= & \sum_{t=0}^{p-1} \sum_{x,y (p)}
\js{t(x^3-x)+1} \js{t(y^3-y)+1} - \sum_{x,y (p)} \js{1}
\nonumber\\ &=& \sum_{x,y (p)} \sum_{t (p)} \js{t(x^3-x)+1}
\js{t(y^3-y)+1} - p^2.
\end{eqnarray}

In Lemma \ref{labquadlegsum} we showed that, if $a$ and $b$ are
not both zero,
\begin{equation}
\twocase{\zsum{t} \js{at^2 + bt + c}\ =\  }{(p-1)\js{a}}{if $p |
b^2 - 4ac$}{-\js{a}}{otherwise.}
\end{equation}

In $A_{2,\mathcal{E}}(p)$ we have
\begin{eqnarray}
a &=& (x^3-x)(y^3-y) = y(x^2-1)x(y^2-1) \nonumber\\ b &=& (x^3-x)
+ (y^3-y) \nonumber\\ c &=& 1 \nonumber\\ \delta(x,y) &=& b^2-4ac
= \Big( (x^3-x) - (y^3-y) \Big)^2.
\end{eqnarray}
We use inclusion / exclusion on $x^3-x$ and $y^3-y$ vanishing.
Assume first that $x^3-x$ equals zero (happens three ways: $x = 0,
\pm 1$). Then we have $\sum_t \js{t(y^3-y)+1}$, which is $3p$ from
our $A_{1,\mathcal{E}}(p)$ computation, giving $3\cdot 3p$.
Similarly we get $3\cdot 3p$ if $y^3-y$ is zero. We subtract the
doubly counted $x^3-x \equiv y^3 - y \equiv 0$ (nine ways), each
of which gives $\sum_t \js{1} = p$. Hence the contribution from at
least one of $x^3-x$ and $y^3-y$ vanishing is $9p$.

Assume $x,y \not\in \{0,\pm 1\}$. When is $\delta(x,y) = (x^3-x) -
(y^3-y) \equiv 0 (p)$?
\begin{eqnarray}
\delta(x,y) &=& (x-y)\cdot(x^2 + xy + y^2 - 1).
\end{eqnarray}
Therefore
\begin{eqnarray}
A_{2,\mathcal{E}}(p) = \sum_{x,y \neq 0,\pm 1 \atop
\delta(x,y)\equiv 0} p \js{(x^3-x)(y^3-y)} - \sum_{x,y \neq 0,\pm
1} \js{(x^3-x)(y^3-y)} + 9p - p^2. \nonumber\\
\end{eqnarray}

Clearly, $\delta(x,y) \equiv 0 (p)$ if $x=y$, which happens $p-3$
times. If $x=y$ then the second factor is $3x^2-1$, which is
congruent to zero at most twice.

When is $\delta_2(x,y) = x^2 + xy + y^2 - 1 \equiv 0$? By the
Quadratic Formula mod $p$,
\begin{eqnarray}
y = \frac{-x \pm \sqrt{4-3x^2}}{2},
\end{eqnarray}
which reduces to finding when $4-3x^2$ is a square mod $p$. We get
two values of $y$ if it is equivalent to a non-zero square, one
value if it is equivalent to zero, and no values if it is not
equivalent to a square. When solving $\delta_2(x,y) \equiv 0 (p)$,
we make sure such $y \not\in \{0,\pm 1\}$. If $y = 0$, $x = \pm
1$; $y = 1$, $x = 0$ or $-1$; $y = -1$, $x = 0$ or $1$. Therefore,
we don't get an excluded $y$ (and similarly if we reverse the
rolls of $y$ and $x$). Thus the number of solutions to
$\delta_2(x,y) \equiv 0 (p)$ is
\begin{eqnarray}
\sum_{x=2}^{p-2} \Bigg[ 1 + \js{4-3x^2} \Bigg] &=& p-3 +
\sum_{x=2}^{p-2} \js{4-3x^2} \nonumber\\ &=& p-6 + \sum_{x (p)}
\js{4-3x^2}.
\end{eqnarray}

We again use Lemma \ref{labquadlegsum}. The discriminant now is
$0^2 - 4\cdot(-3)\cdot4$. For $p \geq 5$, $p$ does not divide the
discriminant, hence this sum is $-\js{-3}$.

Thus, for $x \neq 0, \pm 1$, the number of solutions with $x^2 +
xy + y^2$ $\equiv 1$ is $p-6-\js{-3}$; the number with $x - y
\equiv 0$ is $p-3$. At most two of the pairs $(x,y)$ satisfying
$x^2+xy+y^2-1 \equiv 0 (p)$ also satisfy $x=y$. These pairs
satisfy $3x^2 \equiv 1$, thus, if $\js{3} = 1$ we have doubly
counted two solutions; if it is $-1$, there was no double
counting. Thus, the number of doubly counted pairs is $1 +
\js{3}$, and the total number of pairs is
\begin{eqnarray}
2p - 10 - \js{-3} - \js{3}.
\end{eqnarray}
When $x = y \neq 0, \pm 1$, clearly $\js{(x^3-x)(y^3-y)} = 1$.
Hence these terms contribute $1$.

Consider $x \neq y$ and $x^2+xy+y^2-1 \equiv 0$. Thus $x,y \neq 0,
\pm 1$. Then $y^2-1 \equiv -x(x+y)$ and $x^2-1 \equiv -y(x+y)$ and
\begin{eqnarray}
\js{(x^3-x)(y^3-y)} = \js{x(x^2-1)y(y^2-1)} = \js{x^2 y^2
(x+y)^2}.
\end{eqnarray}

As long as $x \neq -y$, this is $1$. If $x=-y$ then we would have
$x^2-x^2+x^2-1 \equiv 0$. This implies $x = \pm 1$, which cannot
happen as $x,y \neq 0, \pm 1$. Therefore all pairs have their
Legendre factor $+1$, and we need only count how many such pairs
there are. We've previously shown this to be $p + O(1)$, therefore

\begin{eqnarray}
A_{2,\mathcal{E}}(p) &=& p\Bigg[2p - 10 - \js{-3} - \js{3}\Bigg] -
\sum_{x,y \neq 0,\pm 1} \js{(x^3-x)(y^3-y)} + 9p - p^2 \nonumber\\
&=& p^2 - p - \Bigg[\sum_{x (p)} \js{x^3-x}\Bigg]^2 - \js{-3} -
\js{3}.
\end{eqnarray}

As $x^3 - x$ is a non-singular elliptic curve, by Hasse its sum
above is bounded by $4p$. It has complex multiplication and
analytic rank $0$. For $p \equiv 3$ mod $4$ its $a_E(p) = 0$
(change variables $x \to -x$); for the remaining $p$, the angles
of $\frac{a_E(p)}{2\sqrt{p}}$ are uniformly distributed. Hence
$A_{2,\mathcal{E}}(p) = p^2 + O(p)$. \end{proof}

\begin{rek} The reason this calculation succeeds is we have a very
tractable expression for $x(x^2-1)y(y^2-1)$ when $x^2+xy+y^2-1
\equiv 0$ mod $p$. It was non-trivial to find a family with high
rank over $\Q(T)$ and $A_{2,\mathcal{E}}(p)$ computable.\end{rek}

\section{Quadratics Sums of Legendre Symbols}\label{sec:legsums}

For completeness we include proofs of standard sums of Legendre
symbols.

\begin{lem}\label{lablegsum} For $p > 2$
\begin{equation}
S(n)\ =\ \zsum{x} \js{n_1 + x} \js{n_2 + x}\ =\ \Bigg\{ {\ \ p-1 \
\ \ \ \ {\rm if} \ p \ | \ n_1 - n_2 \atop -1 \ \ \ {\rm
otherwise.}}
\end{equation} \end{lem}

\begin{proof} Shifting $x$ by $-n_2$, we need only prove the lemma when
$n_2 = 0$. Assume $(n,p) = 1$ as otherwise the result is trivial.
For $(a,p) = 1$ we have
\begin{eqnarray}
S(n) &\ =\ & \sum_{x=0}^{p-1} \js{n + x} \js{x} \nonumber\\
     & = & \zsum{x} \js{n + a^{-1} x} \js{a^{-1} x} \nonumber\\
     & = & \zsum{x} \js{an + x} \js{x}\ =\ S(an)
\end{eqnarray}
Hence
\begin{eqnarray}
S(n) &\ =\ & \frac{1}{p-1} \osum{a} \zsum{x} \js{an+x} \js{x} \nonumber\\
     & = & \frac{1}{p-1} \zsum{a} \zsum{x} \js{an+x} \js{x} -
     \frac{1}{p-1} \zsum{x} \js{x}^2 \nonumber\\
     & = & \frac{1}{p-1} \zsum{x} \js{x} \zsum{a} \js{an+x} - 1
     \nonumber\\
     & = & 0 - 1 = -1
\end{eqnarray}

Where do we use $p > 2$? We used $\sum_{a=0}^{p-1} \js{an+x} = 0$
for $(n,p) = 1$. This is true for all odd primes (as there are
$\frac{p-1}{2}$ quadratic residues, $\frac{p-1}{2}$ non-residues,
and $0$); for $p=2$, there is one quadratic residue, no
non-residues, and $0$.\end{proof}

\begin{lem}[Quadratic Legendre Sums]\label{labquadlegsum} Assume
$a$ and $b$ are not both zero mod $p$ and $p > 2$. Then
\begin{equation}
\twocase{\zsum{t} \js{at^2 + bt + c}\ =\ }{(p-1)\js{a}}{if $p |
b^2 - 4ac$}{-\js{a}}{otherwise.}
\end{equation}
\end{lem}

\begin{proof} Assume $a \not\equiv 0 (p)$ as otherwise the proof is
trivial. Let $\delta = 4^{-1}(b^2 - 4ac)$. Then
\begin{eqnarray}
\zsum{t} \js{at^2 + bt + c} & = & \zsum{t} \js{a^{-1}} \js{a^2 t^2
+ bat + ac} \\ & = & \zsum{t} \js{a} \js{t^2 + bt + ac}
\nonumber\\ & = & \zsum{t} \js{a} \js{t^2 + bt + 4^{-1}b^2 + ac -
4^{-1} b^2} \nonumber\\ & = & \zsum{t} \js{a} \js{ (t + 2^{-1}b)^2
- 4^{-1}(b^2-4ac)} \nonumber\\ & = & \zsum{t} \js{a} \js{t^2 -
\delta} \nonumber\\ & = & \js{a} \zsum{t} \js{t^2 -\delta}
\nonumber\
\end{eqnarray}
If $\delta \equiv 0 \bmod p$ we get $p-1$. If $\delta = \eta^2,
\eta \neq 0$, then by the Lemma \ref{lablegsum}
\begin{equation}
\zsum{t} \js{t^2 - \delta}\ =\ \zsum{t} \js{t-\eta} \js{t+\eta}\
=\ -1.
\end{equation}
We note that $\zsum{t} \js{t^2 - \delta}$ is the same for all
non-square $\delta$'s (let $g$ be a generator of the
multiplicative group, $\delta = g^{2k+1}$, change variables by $t
\to g^k t$). Denote this sum by $S$, the set of non-zero squares
by $\mathcal{R}$, and the non-squares by $\mathcal{N}$. Since
$\zsum{\delta} \js{t^2 - \delta} = 0$ we have
\begin{eqnarray}
\zsum{\delta} \zsum{t} \js{t^2 - \delta} &\ =\ & \zsum{t} \js{t^2}
+ \sum_{\delta \in \mathcal{R}} \zsum{t} \js{t^2 - \delta} +
\sum_{\delta \in \mathcal{N}} \zsum{t} \js{t^2 - \delta}
\nonumber\\ & = & (p-1) + \frac{p-1}{2}(-1) + \frac{p-1}{2}S\ =\ 0
\end{eqnarray}
Hence $S = -1$, proving the lemma. \end{proof}

\end{document}